\documentclass[12pt]{article}

\usepackage{amsmath,amsthm}
\usepackage{amssymb}
\usepackage{color}
\usepackage{xspace}
\usepackage[colorlinks=true,
linkcolor=green,
filecolor=brown,
citecolor=green]{hyperref}

\def\blue{\textcolor{blue} }
\def\green{\textcolor{green} }
\def\black{\textcolor{black} }

\def\v{\vert}

\def\a{\ensuremath{\mathcal A}\xspace}
\def\bo{\ensuremath{\mathcal B}\xspace}

\def\gf{generating function\xspace}

\renewcommand{\c}[2]{C^{\!\!\!\!\textrm{ \raisebox{1mm}{ {\scriptsize $(#2)$} } }}_{\,#1} \hspace*{-2mm}}

\renewcommand{\b}[2]{\mathcal{B}^{\!\!\!\!\textrm{ \raisebox{1mm}{ {\scriptsize
$(#2)$} \hspace*{-2.8mm}} 
}}_{\,#1}}

\renewcommand{\d}[2]{\mathcal{D}^{\!\!\!\!\textrm{ \raisebox{1mm}{ {\scriptsize $(#2)$} \hspace*{-2.8mm}} }}_{\,#1}}

\newcommand{\f}[2]{\mathcal{F}^{\!\!\!\!\textrm{ \raisebox{1mm}{ {\scriptsize
$(#2)$} \hspace*{-2.8mm}} 
}}_{\,#1}}

\catcode`\@=11

\thicklines
\newskip\Einheit \Einheit=.6cm
\newcount\xcoord \newcount\ycoord
\newdimen\xdim \newdimen\ydim \newdimen\PfadD@cke \newdimen\Pfadd@cke
\def\PfadDicke#1{\PfadD@cke#1 \divide\PfadD@cke by2 
\Pfadd@cke\PfadD@cke \multiply\PfadD@cke by2}
\long\def\LOOP#1\REPEAT{\def\BODY{#1}\ITERATE}
\def\ITERATE{\BODY \let\next\ITERATE \else\let\next\relax\fi \next}
\let\REPEAT=\fi
\def\Punkt{\hbox{\raise-2pt\hbox to0pt{\hss\scriptsize$\bullet$\hss}}}

\def\DuennPunkt(#1,#2){\unskip
  \raise#2 \Einheit\hbox to0pt{\hskip#1 \Einheit
          \raise-1.5pt\hbox to0pt{\hss\tiny$\bullet$\hss}\hss}}
		  
\def\NormalPunkt(#1,#2){\unskip
  \raise#2 \Einheit\hbox to0pt{\hskip#1 \Einheit
          \raise-3pt\hbox to0pt{\hss\large$\bullet$\hss}\hss}}
\def\DickPunkt(#1,#2){\unskip
  \raise#2 \Einheit\hbox to0pt{\hskip#1 \Einheit
          \raise-4pt\hbox to0pt{\hss\Large$\bullet$\hss}\hss}}
\def\Kreis(#1,#2){\unskip
  \raise#2 \Einheit\hbox to0pt{\hskip#1 \Einheit
          \raise-4pt\hbox to0pt{\hss\Large$\circ$\hss}\hss}}
\def\Diagonale(#1,#2)#3{\unskip\leavevmode
  \xcoord#1\relax \ycoord#2\relax
      \raise\ycoord \Einheit\hbox to0pt{\hskip\xcoord \Einheit
         \unitlength\Einheit
         \line(1,1){#3}\hss}}
\def\AntiDiagonale(#1,#2)#3{\unskip\leavevmode
  \xcoord#1\relax \ycoord#2\relax \advance\xcoord by -0.05\relax
      \raise\ycoord \Einheit\hbox to0pt{\hskip\xcoord \Einheit
         \unitlength\Einheit
         \line(1,-1){#3}\hss}}
\def\Pfad(#1,#2),#3\endPfad{\unskip\leavevmode
  \xcoord#1 \ycoord#2 \thicklines\ZeichnePfad#3\endPfad\thinlines}
\def\ZeichnePfad#1{\ifx#1\endPfad\let\next\relax
  \else\let\next\ZeichnePfad
    \ifnum#1=1
      \raise\ycoord \Einheit\hbox to0pt{\hskip\xcoord \Einheit
         \vrule height\Pfadd@cke width1 \Einheit depth\Pfadd@cke\hss}%
      \advance\xcoord by 1
     \else\ifnum#1=2
      \raise\ycoord \Einheit\hbox to0pt{\hskip\xcoord \Einheit
         \unitlength\Einheit
         \line(0,1){1}\hss}
      \advance\xcoord by 0
      \advance\ycoord by 1
 \else\ifnum#1=3
      \raise\ycoord \Einheit\hbox to0pt{\hskip\xcoord \Einheit
         \unitlength\Einheit
         \line(1,1){1}\hss}
      \advance\xcoord by 1
      \advance\ycoord by 1
    \else\ifnum#1=4
      \raise\ycoord \Einheit\hbox to0pt{\hskip\xcoord \Einheit
         \unitlength\Einheit
         \line(1,-1){1}\hss}
      \advance\xcoord by 1
      \advance\ycoord by -1
   \else\ifnum#1=5
      \raise\ycoord \Einheit\hbox to0pt{\hskip\xcoord \Einheit
         \unitlength\Einheit
         \line(2,1){2}\hss}
      \advance\xcoord by 2
      \advance\ycoord by 1
	  \else\ifnum#1=6
      \raise\ycoord \Einheit\hbox to0pt{\hskip\xcoord \Einheit
         \unitlength\Einheit
         \line(2,-1){2}\hss}
      \advance\xcoord by 2
      \advance\ycoord by -1
	  \else\ifnum#1=7
      \raise\ycoord \Einheit\hbox to0pt{\hskip\xcoord \Einheit
         \unitlength\Einheit
         \line(3,1){3}\hss}
      \advance\xcoord by 3
      \advance\ycoord by 1
	  \else\ifnum#1=8
      \raise\ycoord \Einheit\hbox to0pt{\hskip\xcoord \Einheit
         \unitlength\Einheit
         \line(3,-1){3}\hss}
      \advance\xcoord by 3
      \advance\ycoord by -1
    \fi\fi\fi\fi\fi\fi\fi\fi
  \fi\next}
\def\hSSchritt{\leavevmode\raise-.4pt\hbox 
to0pt{\hss.\hss}\hskip.2\Einheit
  \raise-.4pt\hbox to0pt{\hss.\hss}\hskip.2\Einheit
  \raise-.4pt\hbox to0pt{\hss.\hss}\hskip.2\Einheit
  \raise-.4pt\hbox to0pt{\hss.\hss}\hskip.2\Einheit
  \raise-.4pt\hbox to0pt{\hss.\hss}\hskip.2\Einheit}
\def\vSSchritt{\vbox{\baselineskip.2\Einheit\lineskiplimit0pt
\hbox{.}\hbox{.}\hbox{.}\hbox{.}\hbox{.}}}
\def\DSSchritt{\leavevmode\raise-.4pt\hbox to0pt{%
  \hbox to0pt{\hss.\hss}\hskip.2\Einheit
  \raise.2\Einheit\hbox to0pt{\hss.\hss}\hskip.2\Einheit
  \raise.4\Einheit\hbox to0pt{\hss.\hss}\hskip.2\Einheit
  \raise.6\Einheit\hbox to0pt{\hss.\hss}\hskip.2\Einheit
  \raise.8\Einheit\hbox to0pt{\hss.\hss}\hss}}
\def\dSSchritt{\leavevmode\raise-.4pt\hbox to0pt{%
  \hbox to0pt{\hss.\hss}\hskip.2\Einheit
  \raise-.2\Einheit\hbox to0pt{\hss.\hss}\hskip.2\Einheit
  \raise-.4\Einheit\hbox to0pt{\hss.\hss}\hskip.2\Einheit
  \raise-.6\Einheit\hbox to0pt{\hss.\hss}\hskip.2\Einheit
  \raise-.8\Einheit\hbox to0pt{\hss.\hss}\hss}}
\def\SPfad(#1,#2),#3\endSPfad{\unskip\leavevmode
  \xcoord#1 \ycoord#2 \ZeichneSPfad#3\endSPfad}
\def\ZeichneSPfad#1{\ifx#1\endSPfad\let\next\relax
  \else\let\next\ZeichneSPfad
    \ifnum#1=1
      \raise\ycoord \Einheit\hbox to0pt{\hskip\xcoord \Einheit
         \hSSchritt\hss}%
      \advance\xcoord by 1
    \else\ifnum#1=2
      \raise\ycoord \Einheit\hbox to0pt{\hskip\xcoord \Einheit
        \hbox{\hskip-2pt \vSSchritt}\hss}%
      \advance\ycoord by 1
    \else\ifnum#1=3
      \raise\ycoord \Einheit\hbox to0pt{\hskip\xcoord \Einheit
         \DSSchritt\hss}
      \advance\xcoord by 1
      \advance\ycoord by 1
    \else\ifnum#1=4
      \raise\ycoord \Einheit\hbox to0pt{\hskip\xcoord \Einheit
         \dSSchritt\hss}
      \advance\xcoord by 1
      \advance\ycoord by -1
    \fi\fi\fi\fi
  \fi\next}
\def\Koordinatenachsen(#1,#2){\unskip
 \hbox to0pt{\hskip-.5pt\vrule height#2 \Einheit width.5pt depth1 
\Einheit}%
 \hbox to0pt{\hskip-1 \Einheit \xcoord#1 \advance\xcoord by1
    \vrule height0.25pt width\xcoord \Einheit depth0.25pt\hss}}
\def\Koordinatenachsen(#1,#2)(#3,#4){\unskip
 \hbox to0pt{\hskip-.5pt \ycoord-#4 \advance\ycoord by1
    \vrule height#2 \Einheit width.5pt depth\ycoord \Einheit}%
 \hbox to0pt{\hskip-1 \Einheit \hskip#3\Einheit 
    \xcoord#1 \advance\xcoord by1 \advance\xcoord by-#3 
    \vrule height0.25pt width\xcoord \Einheit depth0.25pt\hss}}
\def\Gitter(#1,#2){\unskip \xcoord0 \ycoord0 \leavevmode
  \LOOP\ifnum\ycoord<#2
    \loop\ifnum\xcoord<#1
      \raise\ycoord \Einheit\hbox to0pt{\hskip\xcoord 
\Einheit\Punkt\hss}%
      \advance\xcoord by1
    \repeat
    \xcoord0
    \advance\ycoord by1
  \REPEAT}
\def\Gitter(#1,#2)(#3,#4){\unskip \xcoord#3 \ycoord#4 \leavevmode
  \LOOP\ifnum\ycoord<#2
    \loop\ifnum\xcoord<#1
      \raise\ycoord \Einheit\hbox to0pt{\hskip\xcoord 
\Einheit\Punkt\hss}%
      \advance\xcoord by1
    \repeat
    \xcoord#3
    \advance\ycoord by1
  \REPEAT}
\def\Label#1#2(#3,#4){\unskip \xdim#3 \Einheit \ydim#4 \Einheit
  \def\lo{\advance\xdim by-.5 \Einheit \advance\ydim by.5 \Einheit}%
  \def\llo{\advance\xdim by-.25cm \advance\ydim by.5 \Einheit}%
  \def\loo{\advance\xdim by-.5 \Einheit \advance\ydim by.25cm}%
  \def\o{\advance\ydim by.25cm}%
  \def\ro{\advance\xdim by.5 \Einheit \advance\ydim by.5 \Einheit}%
  \def\rro{\advance\xdim by.25cm \advance\ydim by.5 \Einheit}%
  \def\roo{\advance\xdim by.5 \Einheit \advance\ydim by.25cm}%
  \def\l{\advance\xdim by-.30cm}%
  \def\r{\advance\xdim by.30cm}%
  \def\lu{\advance\xdim by-.5 \Einheit \advance\ydim by-.6 \Einheit}%
  \def\llu{\advance\xdim by-.25cm \advance\ydim by-.6 \Einheit}%
  \def\luu{\advance\xdim by-.5 \Einheit \advance\ydim by-.30cm}%
  \def\u{\advance\ydim by-.30cm}%
  \def\ru{\advance\xdim by.5 \Einheit \advance\ydim by-.6 \Einheit}%
  \def\rru{\advance\xdim by.25cm \advance\ydim by-.6 \Einheit}%
  \def\ruu{\advance\xdim by.5 \Einheit \advance\ydim by-.30cm}%
  #1\raise\ydim\hbox to0pt{\hskip\xdim
     \vbox to0pt{\vss\hbox to0pt{\hss$#2$\hss}\vss}\hss}%
}
\catcode`\@=12

\begin{document}
\newtheorem{theorem}{Theorem}
\newtheorem{defn}[theorem]{Definition}
\newtheorem{lemma}[theorem]{Lemma}
\newtheorem{prop}[theorem]{Proposition}
\newtheorem{cor}[theorem]{Corollary}
\begin{center}
{\Large
A combinatorial interpretation of the Catalan transform of the Catalan numbers    \\ 
}

\vspace*{5mm}

DAVID CALLAN  \\
Department of Statistics  \\
\vspace*{-2mm}
University of Wisconsin-Madison  \\
\vspace*{-2mm}
1300 University Ave  \\
\vspace*{-2mm}
Madison, WI \ 53706-1532  \\
{\bf callan@stat.wisc.edu}  \\
\vspace*{5mm}
\end{center}

\begin{abstract}
The Catalan transform of a sequence $(a_{n})_{n\ge 0}$ is the sequence 
$(b_{n})_{n\ge 0}$ with 
$b_{n}=\sum_{k=0}^{n}\frac{k}{2n-k}\binom{2n-k}{n-k}a_{k}$. 
Here we show that the Catalan transform of 
the Catalan numbers has a simple interpretation: it counts functions 
$f:[1,n]\rightarrow [1,n]$ satisfying the 
condition that, for all $i<j,\ f(j)-(j-i)$ 
is not in the interval $[1,f(i)-1]$.
\end{abstract}

\section{Introduction} \vspace*{-4mm}
The Catalan transform \cite{barry} of a sequence $(a_{n})_{n\ge 0}$ is the sequence 
$(b_{n})_{n\ge 0}$ with 
$b_{n}=\sum_{k=0}^{n}\frac{k}{2n-k}\binom{2n-k}{n-k}a_{k}$,
where $\frac{k}{2n-k}\binom{2n-k}{n-k}$ is interpreted as 1 if $n=k=0$.
The Catalan number is $C_{n}=\frac{1}{n+1}\binom{2n}{n}$.
The purpose of this note is to establish 
the following  combinatorial interpretation.
\begin{theorem}\label{thm:1}
The Catalan transform $(b_{n})_{n\ge 0}$ of 
the Catalan numbers $(C_{n})_{n\ge 0}$ counts functions 
$f:[1,n]\rightarrow [1,n]$ satisfying the 
condition that, for all $i<j,\ f(j)-(j-i)$ 
is not in the interval $[1,f(i)-1]$.   
\end{theorem}

We will represent functions $f:[1,n]\rightarrow [1,n]$ as sequences 
$(u_{i})_{i=1}^{n}$ where $u_{i}=f(i) $ and, hence, $1\le u_{i} \le n 
$ for all $i$.

In Section 2 we review the Catalan numbers, and
in Section 3 we establish some preliminary results about sequences 
$(u_{i})_{i=1}^{n}$  of positive integers that satisfy the key condition
\begin{equation}
\hspace*{40mm}   \textrm{\blue{\framebox{ \black{$\vphantom{A^{A^{A}}_{A_{A}} } u_{j}-(j-i) \notin [1,u_{i}-1]$} } }} \quad \quad \quad 1\le i<j \le n
    \label{eq:1}
\end{equation}
of Theorem \ref{thm:1}. In Section 4 we prove the main result and 
Section 5 presents an extension. 

\section{Generalized Catalan numbers}\label{gencat} \vspace*{-4mm}

The generalized Catalan number $\c{n}{k}$ is defined by $\c{n}{k}=\frac{k+1}{2n+k+1}\binom{2n+k+1}{n}$
with $\c{n}{-1}:=1 $ if $n=0$ and $:=0$ if $n\ge 1$. 
Taking $k=0$ gives the ordinary Catalan number 
$C_{n}:=\c{n}{0}=\frac{1}{2n+1}\binom{2n+1}{n}=\frac{1}{n+1}\binom{2n}{n}$ 
= $\binom{2n}{n}-\binom{2n}{n-1}$. It is well known that $C_{n}$ 
counts Dyck $n$-paths \cite{ec2}.  A Dyck $n$-path has $n$ 
upsteps and $n$ downsteps and its semilength is $n$. 
It is also known \cite{woan} that the sequence $\big(\c{n}{k}\big)_{n\ge 0}$ is the $k$-fold 
convolution of $(C_{n})_{n\ge 0}$ \cite{woan}. It follows, as is well 
known, that 
$\c{n}{k}$ is the number of Dyck $(k+n)$-paths that start with at 
least $k$ upsteps---discard the first $k$ upsteps and then consider the 
decomposition into Dyck paths induced by the last upstep at level $i,\ i=k,k-1,\ldots,1$.

We will use an equivalent formulation of the Catalan transform 
obtained by 
reversing the order of summation in the expression for $b_{n}$. Thus 
\[
b_{n}=\sum_{k=0}^{n}C_{k}^{(n-1-k)}a_{n-k}.
\]

\section{Preliminary results} \vspace*{-4mm}
\begin{prop} \label{bseq}
    Fix a nonnegative integer $k$. Let 
    $\b{n}{k}$ denote the set of sequences $(v_{1},\ldots,v_{n})$ 
    satisfying $1\le v_{i} \le i+k$ and condition $($\ref{eq:1}\,$)$. Then 
    $\v\,\b{n}{k}\,\v = \c{n}{k}$.
\end{prop}
Proof. As noted in Section \ref{gencat},  $\c{n}{k}$ is the number of Dyck $(k+n)$-paths 
that start with at least $k$ upsteps. Given such a path, define 
$v_{i}=1\: + $ semilength of the longest 
Dyck subpath immediately preceding the $(i+k)$-th upstep, $1\le i \le 
n$. This is a bijection to $\b{n}{k}$. 
For example, with $k=2$ and $n=6$, the mandatory $k$ upsteps in the Dyck 
path shown are in blue and the $(i+k)$-th 
upstep is labeled with the corresponding $v_{i}$.

\Einheit=0.6cm
\[
\blue{\Pfad(-9,-2),33\endPfad}
\Pfad(-7,0),3\endPfad
\green{\Pfad(-6,1),34343344\endPfad}
\Pfad(2,1),3444434\endPfad
\SPfad(-9,-2),111111111111111111\endSPfad
\blue{\DuennPunkt(-9,-2)}
\blue{\DuennPunkt(-8,-1)}
\blue{\DuennPunkt(-7,0)}
\DuennPunkt(-6,1)
\DuennPunkt(-5,2)
\DuennPunkt(-4,1)
\DuennPunkt(-3,2)
\DuennPunkt(-2,1)
\DuennPunkt(-1,2)
\DuennPunkt(0,3)
\DuennPunkt(1,2)
\DuennPunkt(2,1)
\DuennPunkt(3,2)
\DuennPunkt(4,1)
\DuennPunkt(5,0)
\DuennPunkt(6,-1)
\DuennPunkt(7,-2)
\DuennPunkt(8,-1)
\DuennPunkt(9,-2)
\Label\o{ \textrm{{\footnotesize 1}}}(-5.7,1.4)
\Label\o{ \textrm{{\footnotesize 2}}}(-3.7,1.4)
\Label\o{ \textrm{{\footnotesize 3}}}(-1.7,1.4)
\Label\o{ \textrm{{\footnotesize 1}}}(-0.7,2.4)
\Label\o{ \textrm{{\footnotesize 5}}}(2.3,1.4)
\Label\o{ \textrm{{\footnotesize 8}}}(7.3,-1.6)
\]

The key observation to show this bijection works is that if an upstep in a Dyck path is immediately preceded by a nonempty Dyck subpath and the steps of the maximal such subpath are colored green, as for the $(6+k)$-th upstep in the example, then the Dyck subpath associated with each green upstep is a subpath of the green path.
\qed

In particular, for $k=0$ we have 
\begin{prop}\label{bcor}
    Let 
    $\bo_{n}=\b{n}{0}$ denote the set of sequences $(v_{1},\ldots,v_{n})$ 
    satisfying $1\le v_{i} \le i$ and condition $($\ref{eq:1}\,$)$. Then 
    $\v\,\bo_{n}\,\v = C_{n}$. \qed
\end{prop}
\textbf{Remark}. The map ``reverse and decrement each entry by 1'' is a 
bijection from $\bo_{n}$ to the inversion codes for 231-avoiding 
permutations of $[n]$. (The inversion code for a permutation 
$(p_{1},\ldots,p_{n})$ of $[n]$ is the sequence $(u_{1},\ldots,u_{n})$ 
where $u_{i}$ is the number of $j\in [i+1,n]$ with $p_{i}>p_{j}$.)

\begin{prop}\label{dseq}
    For $s\ge 0$, let 
    $\d{n}{s}$ denote the set of sequences $(u_{1},\ldots,u_{n})$ 
    satisfying $u_{i}\ge 1$ for all $i$, $u_{n}=n+s$, and condition $($\ref{eq:1}\,$)$. Then 
    $\v\,\d{n}{s}\,\v = \c{n-1}{s}$\ \,.
\end{prop}
Proof. Suppose  $(u_{1},\ldots,u_{n}) \in \d{n}{s}$. Take $i \in [n-1]$. 
If $u_{i} \ge 2$, condition (\ref{eq:1}) 
applied with $j=n$ implies $s+i \le 0$ or $s+i \ge u_{i}$. The first 
inequality cannot hold, so $u_{i} \le s+i$.  But the latter inequality obviously 
also holds when $u_{i}=1$. 
Hence, $(u_{1},\ldots,u_{n-1}) \in \b{n-1}{s}$. Apply Proposition 
\ref{bseq}. \qed

\begin{prop}\label{fseq}
    For $k\ge 0$, let 
    $\f{n}{k}$ denote the set of sequences $(u_{1},\ldots,u_{n})$ 
    satisfying $1\le u_{i} \le n+k+1$ for all $i$, $u_{n}>n$, and condition $($\ref{eq:1}\,$)$. Then 
    $\v\,\f{n}{k}\,\v = \c{n}{k}$.
\end{prop}
Proof. Clearly, $\{(u_{1},\ldots,u_{n}) \in \f{n}{k}\,:\ u_{n}=n+s\} 
\subseteq \d{n}{s}$ for $1\le s \le k+1.$ On the other hand, for 
$(u_{1},\ldots,u_{n}) \in \d{n}{s}$, condition (\ref{eq:1}) implies 
$u_{i} \le s+i \le k+1+i \le n+k+1$. Hence, the reverse inclusion also 
holds and $\{(u_{1},\ldots,u_{n}) \in \f{n}{k}\,:\ u_{n}=n+s\} 
= \d{n}{s}$. The 
result now follows from Proposition \ref{dseq} and the 
identity 
\[
\sum_{s=1}^{k+1}\c{n-1}{s}\ =\c{n}{k},
\]
which can be established, for example, by counting nonnegative paths of $k+n$ 
upsteps and $n$ downsteps by $s=$ number of steps weakly after the last 
downstep.

\section{Main result}\label{main} \vspace*{-4mm}
Let $\a_{n}$ denote the set of sequences $\mathbf{u}=(u_{1},\ldots,u_{n})$ 
    satisfying $1\le u_{i} \le n$ for all $i$, and condition 
    (\ref{eq:1}). The main result, Theorem \ref{thm:1}, can now be 
    stated as 
\begin{theorem}\label{mainthm}
    $\v\,\a_{n}\,\v = \sum_{k=0}^{n-1}C_{k}^{(n-1-k)} C_{n-k}$ for $n\ge 
    1$.
\end{theorem}
Proof. Define the statistic $X$ on $\mathbf{u} \in \a_{n}$ by $X=$ largest $i\in [n]$ 
such that $u_{i}>i$, with $X=0$ if there is no such $i$. Set $\a_{n,k}=\{ 
\mathbf{u} \in \a_{n}\,:\, X(\mathbf{u}) =k\},\ 0\le k \le n-1$.  
We claim  $\v\,\a_{n,k}\,\v = C_{k}^{(n-1-k)} C_{n-k}$, from which the 
Theorem follows.
To see the claim, we have $\a_{n,0}=\bo_{n}$ and Proposition 
(\ref{bcor}) says $\v\,\bo_{n}\,\v = C_{n}$. So the claim holds for 
$k=0$. Now suppose $k\ge 1$ and $(u_{1},\ldots,u_{n}) \in \a_{n,k}$. 
We have $u_{k}>k$ and, applying condition (\ref{eq:1}) with  $i=k$ 
and $j \in [k+1,n]$, we have
\[
u_{j}-(j-k) \notin [1,u_{k}-1].
\]
Hence, $u_{j}-(j-k) \le 0$ or $u_{j}-(j-k) \ge u_{k}$, that is, 
$u_{j}\le j-k$ or $u_{j}\ge j+u_{k}-k>j$. We can't have $u_{j}>j$, by 
the maximality in  the definition of $X(\mathbf{u})=k$, and so $u_{j}\le 
j-k$ for $j \in [k+1,n]$.
Now $(v_{1},\ldots,v_{n-k}):=(u_{k+1},\ldots,u_{n})$ inherits condition 
(\ref{eq:1}) and $v_{i} \le i$ for $ 1\le i \le n-k$. Thus 
$(u_{k+1},\ldots,u_{n}) \in \bo_{n-k}$.

Also, for $1\le i \le k$, we have $u_{i} \le n =k+(n-k-1)+1$ and 
$u_{k}>k$, in other words, $(u_{1},\ldots,u_{k}) \in \mathcal{F}_{k}^{(n-k-1)}$.
The map $(u_{1},\ldots,u_{n}) \rightarrow 
\{(u_{1},\ldots,u_{k}),(u_{k+1},\ldots,u_{n})\}$ is in fact a 
bijection from $\a_{n,k}$ to the Cartesian product $\mathcal{F}_{k}^{(n-k-1)} 
\times \bo_{n-k}$.  The claim now follows from the counting results 
of Propositions \ref{bcor} and \ref{fseq}. \qed

\section{Extensions} \vspace*{-4mm}
Similar methods establish a combinatorial interpretation for the 
$k$-shifted Catalan sequence $(0,0,\ldots,0,1,1,2,5,14,\ldots) $ with 
$k$ initial 0's.
\begin{theorem}
   Let $(b_{n})_{n\ge 0}$ be the Catalan transform of the $k$-shifted Catalan 
   sequence. Then $b_{n}=0$ for $n \le k-1,\ b_{k}=1$ and, for 
   $n\ge k+1,\ b_{n}$ is the number of sequences 
   $(u_{1},\ldots,u_{n-k})$ 
    satisfying $1 \le u_{i} \le n$ and condition $($\ref{eq:1}\,$)$. 
\end{theorem}
Stefan Forcey \cite{Forcey} gives a combinatorial 
interpretation for the Catalan transform of the $1$-shifted Catalan sequence 
(\htmladdnormallink{A121988}{http://oeis.org/A121988} in The On-Line Encyclopedia 
of Integer Sequences \cite{oeis}):  $b_{n}$ is the number of vertices of the $n$-th multiplihedron.

\end{document}

\begin{abstract}
The Catalan transform of a sequence $(a_{n})_{n\ge 0}$ with \gf $A(x)= \sum_{k\ge 0}a_{n}x^{n}$ has \gf $A\big(xC(x)\big)$, where $C(x)$ is the \gf  for the Catalan numbers.
Here we show that the Catalan transform of 
the Catalan numbers has a simple interpretation: it counts functions 
$f:[1,n]\rightarrow [1,n]$ satisfying the 
condition that, for all $i<j,\ f(j)-(j-i)$ 
is not in the interval $[1,f(i)-1]$.
\end{abstract}

We note in passing that the sum $\sum_{k=0}^{n}C_{k}^{(n-1-k)}$ 
of the coefficients in the transform is $C_{n}$, because this sum counts 
Dyck $n$-paths by the length $n-k$ of the first ascent---``delete first peak'' 
is a bijection to the  Dyck paths counted by $C_{k}^{(n-1-k)}$.